\documentclass[12pt]{article}

\usepackage{CJK,CJKnumb,amsmath}
\usepackage{amsfonts}
\usepackage{amsthm}
\usepackage{geometry}
\usepackage{mathrsfs}
\usepackage{amsmath}
\usepackage{amssymb}
\usepackage{amsfonts}
\usepackage[percent]{overpic}
\usepackage{bm}
\usepackage[all]{xy}
\usepackage{graphicx}
\usepackage{subfigure}
\usepackage{latexsym}
\usepackage[colorlinks, linkcolor=blue, anchorcolor=blue, citecolor=blue]{hyperref}
\usepackage{epigraph}
\usepackage{hyperref}
\usepackage{fancyhdr}
\usepackage{comment}

\usepackage{mathtools}
\mathtoolsset{showonlyrefs}

\numberwithin{equation}{section}
\usepackage{color}

\newtheorem{theorem}{Theorem}[section]

\newtheorem{definition}{Definition}[section]

\newtheorem{lemma}{Lemma}[section]
\theoremstyle{definition}
\newtheorem{remark}{Remark}[section]

\newtheorem*{ack}{Acknowledgement}

\def\re{\operatorname{Re}}
\def\Im{\operatorname{Im}}

\def\dist{\operatorname{dist}}

\def\c{\operatorname{\mathbb C}}

\def\j{\operatorname{\mathcal{J}}}

\def\b{\operatorname{\mathcal{B}}}

\def\p{\operatorname{\mathcal{P}}}

\begin{document}
\title{Ergodic exponential maps with escaping singular behaviours}
\author{Weiwei Cui and Jun Wang}
\date{}
\maketitle

\begin{abstract}
We construct exponential maps for which the singular value tends to infinity under iterates while the maps are ergodic. This is in contrast with a result of Lyubich from 1987 which tells that $e^z$ is not ergodic.


\end{abstract}

\section{Introduction and main results}

In this paper we are concerned with the exponential family
\[
f_{\lambda}(z)=\lambda e^z,\quad \lambda\in\c\setminus\{0\}.
\]
As one of the most well studied families in transcendental dynamics, it has provided many further inspirations for exploring more general entire functions. There have been a lot of works on dynamical properties of exponential functions; see, for instance, \cite{baker2, devaney8, eremenko-lyubich, rempethesis, rempe-schleicher} and references therein. 

In this paper, we study exponential maps for which the singular value $0$ tends to $\infty$ under iterates. We call them \emph{escaping maps} for simplicity, and the corresponding $\lambda$ is called an escaping parameter. For example, $e^z$ (i.e., $\lambda=1$) is such a map. Escaping parameters play a fundamental role in the study of the structure of the parameter space of the exponential family. A classification of such maps can be found in \cite{foerster}. Our intention is instead to study their measurable dynamics. Here we will focus on ergodicity of escaping maps. Recall that an exponential map is ergodic with respect to the Lebesgue measure if any invariant set has either full measure or zero measure. 

Not much is known about ergodicity of escaping exponential maps, except for a result of Lyubich from 1987 \cite{lyubich} saying that $e^z$ is not ergodic, answering a question of Sullivan. A natural question to ask is whether Lyubich's result holds for all escaping maps. This is plausible if the singular value escapes sufficiently fast. However, we show below that this is not always the case.

\begin{theorem}\label{thm1}
There exists an escaping parameter $\lambda$ such that $f_{\lambda}$ is ergodic.
\end{theorem}

The idea of construction is through perturbations in the parameter space of exponential maps. Recall that $f_{\lambda}$ is post-singularly finite if the forward orbit of $0$ is finite. Perturbing any post-singularly finite parameter we find a sequence of post-singularly finite exponential maps which approximate our desired escaping one. Ergodicity of the escaping map is ``inherited'' from these maps.

The singular behaviour of the escaping map constructed in the above theorem can be described as follows: The singular value $0$ first gets close to some repelling fixed point, stays there for many iterates, then gets mapped close to another repelling fixed point which is further right and stays there for even more iterates, and so forth.

\medskip

 The construction actually tells that such ergodic escaping maps can be found near any post-singularly finite parameters. Let $\b$ denote the bifurcation locus of the exponential family; i.e., $\b=\{\lambda\in\c\setminus\{0\}: \{f_{\lambda}^{n}(0)\}~\text{is not normal at}~ \lambda\}$. Since post-singularly finite parameters are dense in $\b$,  the following result follows directly.

\begin{theorem}\label{thm2}
Ergodic escaping maps are dense in the bifurcation locus $\b$.
\end{theorem}

\section{Construction}

In this section we construct an escaping exponential map for which the singular value “slowly” goes to $\infty$ under iterates. To begin with, we shall need some preliminary results and notations concerning perturbations in the parameter space. 

\smallskip

Let $f_{\lambda}$ be an exponential map. The Julia set of $f_{\lambda}$ is denoted by $\j(f_{\lambda})$. For a post-singularly finite map, the singular value $0$ will eventually fall into a repelling cycle. So the Julia set of a post-singularly finite exponential map is the whole plane.  This is also true for escaping exponential maps; see \cite[Corollary 1]{baker2}.

\smallskip

\noindent \emph{Some notations}. Let $\c$ be the complex plane and $\mathbb{Z}$ the set of integers. The real and imaginary parts of $z\in\c$ are denoted by $\re(z)$ and $\Im(z)$, respectively. We shall use $D(z,r)$ for the disk of radius $r$ centered at $z$ and by $\overline{D}(z,r)$ we mean the closure of $D(z,r)$. We will use $f'(x)$ to denote the derivative of a function $f$ at $x$. We also use $x\gtrsim y$ to mean that $x\geq Cy$ for some absolute constant $C>0$.

\subsection{Some lemmas}

Let $\lambda_0$ be a post-singularly finite parameter and $D(\lambda_0, r)$ a parameter disk around $\lambda_0$ of radius $r>0$. We shall consider the following functions
\begin{equation}\label{xi}
\xi_n(\lambda)=f_{\lambda}^{n}(0)
\end{equation}
which trace the singular orbits for all parameters in $D(\lambda_0, r)$ and relate the parameter space to the phase space.

The expansion along the post-singular set of $f_{\lambda_0}$ implies the existence of a holomorphic motion $h: D(\lambda_0, r)\times \p(f_{\lambda_0})\to \c$ such that for all $\lambda\in D(\lambda_0, r)$ one has $h_{\lambda}\circ f_{\lambda_0}=f_{\lambda}\circ h_{\lambda}$ for all $z\in\p(f_{\lambda_0})$. Put $\mu_j(\lambda)=h_{\lambda}(\xi_j(\lambda_0))$. Moreover, the expansion on $\p(f_{\lambda_0})$ also implies the existence of a neighborhood of $\p(f_{\lambda_0})$ such that the function has expansion in this neighborhood. Thus one can take a number $\delta>0$ so that the set of points $z$ with $\dist(z, h_{\lambda}(\p(f_{\lambda_0})))<10\delta$ is contained in this neighborhood for all $\lambda\in D(\lambda_0, r)$.

Let $\mathcal{N}:=D(\p(f_{\lambda_0}), 10\delta)=\bigcup_{w\in \p(f_{\lambda_0})}D(w, 10\delta)$. The following lemma, which was proved in \cite[Lemma 3.6, Lemma 3.7]{aspenberg-cui}, is key to our construction.

\begin{lemma}\label{dislarge}
Let $\lambda_0$ be a post-singularly finite parameter. For any $\varepsilon>0$ there exist $\delta>0$ and $r>0$ sufficiently small such that for any $\lambda_1,\lambda_2\in D(\lambda_0, r)$, if $|\xi_j(\lambda_i)-\mu_j(\lambda_i)|\leq \delta$ for $i=1,2$ and for $j\leq n$, then
\begin{equation}\label{distort}
\left|\frac{\xi'_n(\lambda_1)}{\xi'_n(\lambda_2)}-1 \right|<\varepsilon.
\end{equation}
Moreover, there exists $S>0$ such that for all $r>0$ small enough one can find $N_0\in\mathbb{N}$ such that $D(\xi_{N_0}(\lambda_0), S/4)\subset \xi_{N_0}(D(\lambda_0, r))\subset \mathcal{N}$ and $\xi_{N_0}(D(\lambda_0, r))$ has diameter at least $S$.
\end{lemma}

We say that the parameter disk $D(\lambda, r)$ has reached \emph{large scale} $S$ at time $n$ if there exists $n$ so that $\xi_{n}(D(\lambda, r))$ has diameter at least $S$.

\smallskip

The following result, which is often called the \emph{blowing up property of Julia sets}, can be found in \cite[Lemma 2.2]{baker15}; see also \cite{bergweiler1}.

\begin{lemma}\label{blowup}
Let $f_{\lambda}$ be an exponential map. Let $U\subset\c$ be a neighborhood of $z\in\j(f_{\lambda})$ and let $K\subset\c$ be compact not containing $0$. Then there exists an integer $N\in\mathbb{N}$ depending on $K$ such that $f_{\lambda}^{n}(U)\supset K$ for all $n\geq N$.
\end{lemma}

For a post-singularly finite map $f_{\lambda}$, the preperiod of $f_{\lambda}$ is denoted by $n_{\lambda}$, i.e., $f_{\lambda}^{n_{\lambda}}(0)$ is periodic under $f_{\lambda}$. A post-singularly finite map $f_{\lambda}$ is said to be \emph{post-singularly fixed}, if $f_{\lambda}^{n_{\lambda}}(0)$ is a fixed point, which we denote by $\eta_{\lambda}$.

\begin{definition}
A post-singularly fixed map $f_{\lambda}$ is said to be special, if there exists $k\in\mathbb{Z}$ so that
\begin{equation}\label{posit}
|\Im(\eta_{\lambda})-\arg(\lambda)-2k\pi|<\pi
\end{equation}
and
\begin{equation}\label{reall}
\re \eta_{\lambda}\gtrsim \max_{n<n_{\lambda}}\re f_{\lambda}^{n}(0).
\end{equation}
\end{definition}

We first show how to find a special post-singularly fixed parameter near post-singularly finite ones.
\begin{lemma}\label{lepre}
Let $\lambda_0$ be a post-singularly finite parameter. Then for any $r>0$ sufficiently small, there exists a special post-singularly fixed parameter $\lambda_1\in D(\lambda_0, r)$.
\end{lemma}

\begin{proof}
Recall that the functions $\xi_n$ are defined in \eqref{xi}. By Lemma \ref{dislarge}, given $S>0$ one can find $N_0\in\mathbb{N}$ depending on $r$ and  $D_0\subset D(\lambda_0, r)$ such that $\xi_{N_0}(D_0)=D(\xi_{N_0}(\lambda_0), S/4)$. 

Let $x_0>4\pi$ be any positive real number. By Lemma \ref{blowup}, there exists a minimal positive integer $N_1$ such that $f_{\lambda_0}^{N_1}(\xi_{N_0}(D_0))$ contains a square $Q_0$ of sidelength $2\pi$ whose center is 
\[z_0=x_0 -i\arg(\lambda_0),\quad\text{where}\quad \arg(\lambda_0)\in [0,2\pi).\]
Put $\mathcal{L}:=\{z:\, \Im(z)=-\arg(\lambda_0)\}$, i.e., the horizontal line passing through the point $-i\arg(\lambda_0)$. Define
$$\hat{x}_0=\max\left\{\max_{j\leq N_0}\re f_{\lambda_0}^{j}(0)+2\pi,\, \sup_{0\leq j\leq N_1}\sup_{z\in \xi_{N_0}(D_0)}\re f_{\lambda_0}^{j}(z)\right\}.$$

\begin{figure}[htbp]
	\centering
	\includegraphics[width=12cm]{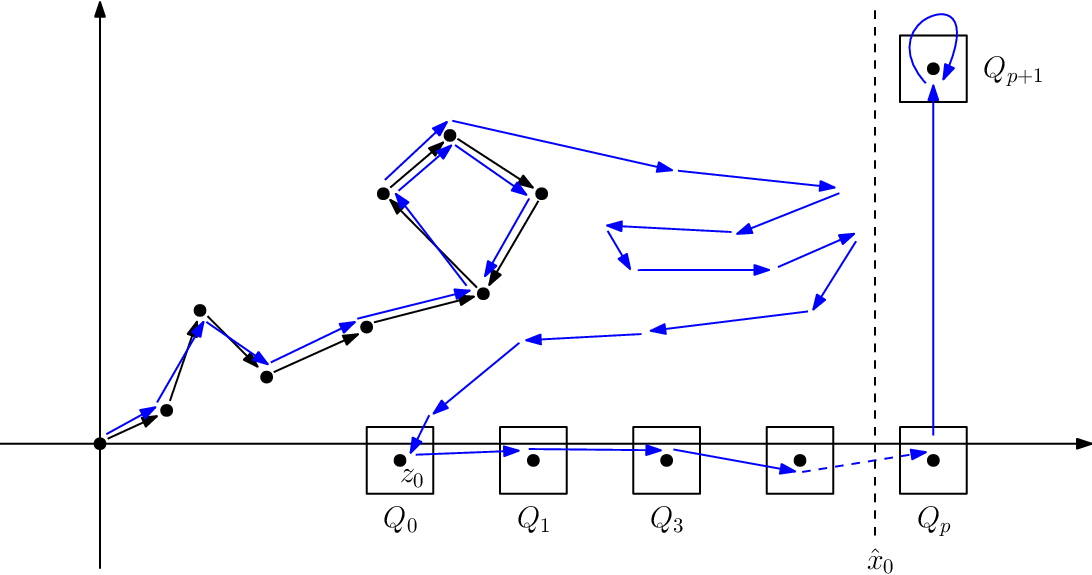}
	\caption{Black ones denotes the singular orbit of the starting map $f_{\lambda_0}$. The blue ones illustrates singular orbit of the desired special post-singularly fixed parameter.}
	\label{pic1}
\end{figure}

Note that $f_{\lambda_0}(Q_0)$ is an annulus $A_0$ of inner radius $|\lambda_0| e^{x_0-\pi}$ and outer radius $|\lambda_0| e^{x_0+\pi}$. So, there exists another square $Q_1\subset A_0$ of sidelength $2\pi$ with center $z_1\in \mathcal{L}$ and $\re(z_1)>0$. Notice again $f_{\lambda_0}(Q_1)$ is also an annulus which we denote it by $A_1$. Then similarly$A_1$ contains a square $Q_2$ of sidelength $2\pi$ centering at a point $z_2\in\mathcal{L}$. By repeating this process, one can find a finite number of squares $Q_j$ centering at $z_j\in\mathcal{L}$ such that $f_{\lambda_0}(Q_j)\supset Q_{j+1}$. By the choice of $x_0$ and $\mathcal{L}$ and by choosing $z_j$ carefully, we could easily achieve that $\re z_{j+1}> e^{\re z_j}$. Now one can see that there exists a number $p\in\mathbb{N}$ so that $\re z_p > \hat{x}_0+2\pi$. Since $f_{\lambda_0}(Q_p)$ is a large annulus, it must contain a $2k\pi i$ translate $Q_{p+1}$ of $Q_p$ for some $k\in\mathbb{Z}$. Note also that $f_{\lambda_0}(Q_{p+1})=f_{\lambda_0}(Q_p)$.

By making perturbations sufficiently small, i.e., taking $r$ sufficiently small, one can achieve that for all $\lambda\in D_0$,
$$f_{\lambda}^{N_1+p+2}(D(\xi_{N_0}(\lambda_0), S/4))\supset Q_{p+1}.$$
Let $V$ be the component of $f_{\lambda}^{-(N_1+p+2)}(Q_{p+1})$ which is contained in $D(\xi_{N_0}(\lambda_0), S/4)$. Also, let $D'_0\subset D_0$ be the set of parameters  which are mapped by $\xi_{N_0}$ into $V$; i.e., $\xi_{N_0}(D'_0)=V$. It follows from the above construction that all parameters $\lambda\in D'_0$ have similar singular behaviours up to $N_0+N_1+p+2$ iterates: They first follow closely the singular orbit of $f_{\lambda_0}$ for $N_0$ iterates, hit the square $Q_0$ by using $N_1$ iterates, and then are mapped into $Q_1,\dots, Q_{p+1}$ accordingly, and then stay in $Q_{p+1}$ for one more iterate. See Figure \ref{pic1} for an illustration of the construction.

We still need to prove the existence of a post-singularly fixed parameter $\lambda_1\in D'_0$. This follows from the above  construction. Let $B$ be the preimage of $\overline{Q}_{p+1}$ under $f_{\lambda}^{-1}$, where the inverse branch is chosen such that $f_{\lambda}(Q_{p+1})\supset Q_{p+1}$. Since $B$ is compactly contained in $Q_{p+1}$ by construction, we see that $f_{\lambda}^{-1}$ must have an attracting fixed point in $B$. In other words, $f_{\lambda}$ has a repelling fixed point in $B$. Pulling back to the parameter space under the function $\xi_{N_0+N_1+p+2}$, this gives a post-singularly fixed parameter $\lambda_1$ which is contained in $D'_0\subset D(\lambda_0, r)$. By the choice of $Q_{p+1}$ we see that \eqref{reall} must be satisfied by this parameter. Therefore, the parameter $\lambda_1$ is special. 
\end{proof}

\begin{lemma}\label{gd}
Let $\lambda_0$ be a special post-singularly fixed parameter and $r>0$ small enough. Then for any $x,\,y$ satisfying $x>y\geq \re\eta_{\lambda_0}$, there exist a special post-singularly fixed parameter $\lambda_1\in D(\lambda_0, r)$, $\tilde{N}$ and $N\in\mathbb{N}$ such that $\tilde{N}<N$ and
\begin{equation}\label{lowerb}
\re f_{\lambda_1}^{n}(0)\,
\begin{cases}
\, \geq x\quad \text{~for~}\quad n\geq N;\\
 \, \leq y\quad \text{~for~}\quad n< \tilde{N}.
\end{cases}
\end{equation}
Moreover, $\tilde{N}$ can be taken arbitrarily large by choosing sufficiently small $r$.
\end{lemma}

The choice of $y$ and thus $\tilde{N}$ make sure that the singular orbit of $\lambda_1$ does not go back to some left half-plane.

\begin{proof}
The idea of proof is similar to that of Lemma \ref{lepre}, except that Lemma \ref{blowup} is not used here. Instead, we shall use property \eqref{posit} of special post-singularly fixed parameters. In this way, we will gain control of the whole singular orbit of the desired parameter, while we usually do not get much information from Lemma \ref{blowup}. We give a sketch of proof for completeness.

As before, Lemma \ref{dislarge} ensures the existence of $N_0$ such that $\xi_n(D(\lambda_0, r))$ reaches the large scale $S$ with bounded distortion. Since $f_{\lambda_0}$ is post-singularly fixed, the large scale must happen around $\eta_{\lambda_0}$. Now that $\lambda_0$ is special, \eqref{posit} holds. Continuing iteration, one would see that there exists $N_1$ so that $\xi_{N_0+N_1}(D(\lambda_0, r))\supset Q_0$, where $Q_0$ is a square of sidelength $2\pi$ around $\eta_{\lambda_0}$. Here $Q_0$ plays the same role as in the proof of previous lemma.

For any $x>y\geq \re \eta_{\lambda_0}$, one can find a finite number, say $p+q$, squares $Q_j$ of sidelength $2\pi$ centering on the line $\{z: \Im(z)=-\arg(\lambda_0)\}$ such that $f_{\lambda_0}(Q_j)\supset Q_{j+1}$ and $Q_j$ for $j\leq p-1$ belongs to the left half-plane $\{z: \re(z)\leq y\}$, $Q_j$ for $p\leq j\leq p+q-1$ belongs to the vertical strip $\{z: y<\re(z)<x\}$, while $Q_{p+q}$ is the only square in the right half-plane $\{z: \re(z)\geq x\}$. Since $f_{\lambda_0}(Q_{p+q})$ is a large annulus, it also contains some $2k\pi i$ translates of $Q_{p+q}$. The rest of construction is the same as in the last lemma. So we omit details.

When taking $r$ sufficiently small, it will take longer iterates to reach the large scale. Thus $\tilde{N}$ will also be large.
\end{proof}

\begin{remark}
It is not difficult to see from the proof of the above lemma that the singular orbit of $\lambda_1$ in the vertical strip $\{z: y<\re(z)<x\}$ are asymptotic to some horizontal line. This is not necessary. In principle, one can make the singular orbit of $\lambda_1$ stays in this strip as many iterates as we want. This is achieved by choosing as many squares as we wish.
\end{remark}

\subsection{The construction and proof of ergodicity}

We start the construction of an ergodic exponential map in this part. Fix $R>0$ large throughout and let $\{\delta_k\}$ be a sequence of positive numbers satisfying $\delta_k\to 0$ as $k\to\infty$. We also fix a strictly increasing sequence $\{x_n\}$ so that $x_n\to\infty$ as $n\to\infty$ and $x_0>0$ sufficiently large.

The following property of post-singularly finite maps is well known.

\begin{lemma}\label{psfergodic}
For any post-singularly finite exponential map, almost every point accumulates everywhere in the whole plane.
\end{lemma}

\begin{proof}
A result of Bock \cite{bock} says that, for any post-singularly finite exponential map almost every point in the plane either
\begin{itemize}
\item[(i)]  tends to $\infty$ under iterates, or
\item[(ii)] accumulates everywhere in the whole plane.
\end{itemize}
However, the set of points tending to $\infty$ under iterates, called the escaping set, has Lebesgue measure zero by a result of Eremenko and Lyubich \cite[Section 7]{eremenko-lyubich}. Therefore,  post-singularly finite maps fall into category $(ii)$.
\end{proof}

Let $\lambda_0$ be a post-singularly finite parameter and $r_0>0$. Then by Lemma \ref{lepre}, there exists a special post-singularly fixed parameter $\lambda_1\in D(\lambda_0, r_0)$. So there exists $\eta_1>0$ such that $\bigcup_{z\in \p(f_{\lambda_1})}D(z,\eta_1)$ has Lebesgue measure less than $\delta_1$. By Lemma \ref{gd} one can find a special post-singularly fixed parameter $\lambda_2\in D(\lambda_1, r_1)$ for a fixed small $r_1>0$. Moreover, there exists $\tilde{N}_2$ and $N_2$ satisfying $\tilde{N}_2<N_2$ so that with $y_2:=(x_1+x_2)/2$ the following holds:
\begin{equation}
\re f_{\lambda_2}^{n}(0)\,
\begin{cases}
\, \geq x_2\quad \text{~for~}\quad n\geq N_2;\\
 \, \leq y_2\quad \text{~for~}\quad n < \tilde{N}_2.
\end{cases}
\end{equation}
Given $\delta_2$, a number $\eta_2$ can now be chosen so that $\bigcup_{z\in \p(f_{\lambda_2})}D(z,\eta_2)$ has Lebesgue measure less than $\delta_2$.

We continue the above process by induction. Suppose now that we have constructed $r_{k-1}$, a special post-singularly fixed parameter $\lambda_k\in D(\lambda_{k-1}, r_{k-1})$, $\tilde{N}_k$ and $N_k$ satisfying $\tilde{N}_k<N_k$ so that with $y_k:=(x_{k-1}+x_k)/2$ we have
\begin{equation}\label{esp}
\re f_{\lambda_k}^{n}(0)\,
\begin{cases}
\, \geq x_k\quad \text{~for~}\quad n\geq N_k;\\
 \, \leq y_k\quad \text{~for~}\quad n< \tilde{N}_k.
\end{cases}
\end{equation}
So we can choose $\eta_k$ such that $\bigcup_{z\in \p(f_{\lambda_k})}D(z,\eta_k)$ has Lebesgue measure less than $\delta_k$.

Fix a point $w\in D(0,R)\setminus \bigcup_{z\in \p(f_{\lambda_k})}D(z,\eta_k)$ and consider the pullbacks of $D(w, \eta_k/4)$ under $f_{\lambda_k}$. By Lemma \ref{psfergodic}, there is $M_k\in\mathbb{N}$ so that 
$$D(0,R)\setminus \bigcup_{j=0}^{M_k} f_{\lambda_k}^{-j}(D(w,\eta_k/4))$$
has Lebesgue measure at most $\delta_k$. The choice of $w$ implies that $f_{\lambda_k}^{-j}$ has bounded distortion on $D(w,\eta_k/2)$. Moreover, $f_{\lambda_k}$ being post-singularly finite implies that the any branch of the inverse function of $f_{\lambda_k}$ is a strict contraction with respect to the hyperbolic metric of $\mathbb{C}\setminus \p(f_{\lambda_k})$. So one can find a number $\varepsilon_k>0$ such that for all $j\leq M_k$, the component of $f_{\lambda_k}^{-j}(D(w,\eta_k/4))$ intersecting with $D(0,R)$ has diameter at most $\varepsilon_k$ for large $k$. Moreover, $\varepsilon_k\to 0$ as $k\to\infty$.

By making perturbations sufficiently small, i.e., choosing $r_k$ small enough, one can achieve that the Lebesgue measure of 
$$D(0,R)\setminus \bigcup_{j=0}^{M_k} f_{\lambda'}^{-j}(D(w,\eta_k/4))$$
is less than $2\delta_k$ for all $\lambda'\in D(\lambda_k, r_k)$, and $f_{\lambda'}^{-j}$ also has bounded distortion (with a possibly larger constant) for $1\leq j\leq M_k$. Moreover, any component of $f_{\lambda'}^{-j}(D(w,\eta_k/4))$ intersecting with $D(0,R)$ has diameter at most $2\varepsilon_k$. Recall that $n_{\lambda_k}$ is the preperiod of $0$ under $f_{\lambda_k}$.  So if $r_k$ is chosen small so that $D(\lambda_k, r_k)\subset D(\lambda_{k-1}, r_{k-1})$, \eqref{esp} holds similarly:
\begin{equation}
\re f_{\lambda'}^{n}(0)\,
\begin{cases}
\, \geq x_k\quad \text{~for~}\quad N_k \leq n\leq n_{\lambda_k};\\
 \, \leq y_k\quad \text{~for~}\quad n< \tilde{N}_k.
\end{cases}
\end{equation}
By Lemma \ref{gd} we can find a special post-singularly fixed parameter $\lambda_{k+1}\in D(\lambda_k, r_k)$ so that \eqref{esp} holds with $k$ replaced by $k+1$. 

Continuing inductively, since $r_k\to 0$ as $k\to\infty$ by construction one obtains a sequence of parameters $\lambda_k$ tending to a limiting parameter $\lambda$, which must be escaping since $x_k$ and hence $y_k$ tend to $\infty$ as $k\to\infty$.

\bigskip

Now we prove that $f_{\lambda}$ is ergodic. Otherwise, there exists an invariant set $A$ such that both $A$ and $B:=\c\setminus A$ have positive Lebesgue measure. Note that $f_{\lambda}$ is escaping, so the Julia set of $f_{\lambda}$ is the whole plane; see, for instance, \cite[Corollary 1]{baker2}. By Lemma \ref{blowup}, both sets $D_{A}:=D(0,R)\cap A$ and $D_B:=D(0,R)\cap B$ have positive measure. 

Since $\bigcup_{z\in \p(f_{\lambda_k})}D(z,\eta_k)$ has Lebesgue measure less than $\delta_k$ which tends to $0$ as $k\to\infty$, for all $k$ large we could take a Lebesgue density point $\alpha\in D_A$ such that 
$$\alpha\not\in \bigcup_{z\in \p(f_{\lambda_k})}D(z,\eta_k).$$

 Therefore, $D(0,R)\setminus \bigcup_{0\leq j\leq M_k} f_{\lambda_k}^{-j}(D(\alpha,\eta_k/4))$ has Lebesgue measure at most $\delta_k$ which tends to $0$ as $k\to\infty$. So, there also exists a Lebesgue density point $\beta\in D_B$ so that for all $k$ large,
$$\beta\in \bigcup_{0\leq j\leq M_k} f_{\lambda_k}^{-j}(D(\alpha,\eta_k/4)).$$
In particular, $\beta$ belongs to some component of $f_{\lambda_k}^{-p_k}(D(\alpha,\eta_k/4))$ for some $p_k\leq M_k$. Since $f_{\lambda}^{-p_k}$ also has bounded distortion on $D(\alpha,\eta_k/2)$ and its diameter $\leq 2\varepsilon_k\to 0$ as $k\to\infty$, one can find a definite portion of $A$ near $\beta$. This contradicts with $\beta$ being a Lebesgue density point of $B$. So, $f_{\lambda}$ is ergodic.

\begin{remark}
As mentioned in the introduction, it is plausible that fast escaping speed of the singular value implies non-ergodicity of the function. Here we constructed an ergodic exponential function with slowly escaping singular value. However, it seems difficult to find a condition on the speed of escaping which implies ergodicity of exponential maps. A natural question is whether there exists a ``slowly'' escaping exponential map which is non-ergodic.
\end{remark}

\medskip

\begin{ack}
The first author was partially supported by the Natural Science Foundation of China (No. 12401105), Shandong Provincial Natural Science Fund for Excellent Young Scientists Program (Overseas)  (No. 2025HWYQ-021), Qingdao Natural Science Foundation (No. 24-4-4-zrjj-8-jch), and also a grant from Vergstiftelsen, Sweden. The second author was partially supported by the Natural Science Foundation of China (No.12471072). We would like to thank Jie Ding and Leticia Pardo-Sim\'on for useful comments. We also thank the referee for a detailed comments and valuable suggestions.
\end{ack}

%\bibliography{references}

\begin{thebibliography}{Rem03}

\bibitem[AC24]{aspenberg-cui}
M.~Aspenberg and W.~Cui, \emph{Perturbations of exponential maps: Non-recurrent
  dynamics}, J. Anal. Math. \textbf{153} (2024), no.~2, 759--775.

\bibitem[Bak84]{baker15}
Irvine~Noel Baker, \emph{Wandering domains in the iteration of entire
  functions}, Proc. London Math. Soc. (3) \textbf{49} (1984), no.~3, 563--576.

\bibitem[Ber93]{bergweiler1}
W.~Bergweiler, \emph{Iteration of meromorphic functions}, Bull. Amer. Math.
  Soc. \textbf{29} (1993), no.~2, 151--188.

\bibitem[Boc96]{bock}
Heinrich Bock, \emph{On the dynamics of entire functions on the {J}ulia set},
  Results Math. \textbf{30} (1996), no.~1-2, 16--20. 

\bibitem[BR84]{baker2}
Irvine~Noel Baker and Philip~J. Rippon, \emph{Iteration on exponential
  functions}, Ann. Acad. Sci. Fenn., Ser. A. I. Math. \textbf{9} (1984),
  49--77.

\bibitem[DK84]{devaney8}
Robert~L. Devaney and Micha{\l} Krych, \emph{Dynamics of {${\rm exp}(z)$}},
  Ergodic Theory Dynam. Systems \textbf{4} (1984), no.~1, 35--52.

\bibitem[EL92]{eremenko-lyubich}
A.~{\`E}. Er\"{e}menko and M.~Yu. Lyubich, \emph{Dynamical properties of some
  classes of entire functions}, Ann. Inst. Fourier (Grenoble) \textbf{42}
  (1992), no.~4, 989--1020. 

\bibitem[FRS08]{foerster}
Markus F\"{o}rster, Lasse Rempe, and Dierk Schleicher, \emph{Classification of
  escaping exponential maps}, Proc. Amer. Math. Soc. \textbf{136} (2008),
  no.~2, 651--663. 

\bibitem[Lyu87]{lyubich}
M.~Yu. Lyubich, \emph{The measurable dynamics of the exponential}, Sibirsk.
  Mat. Zh. \textbf{28} (1987), no.~5, 111--127.

\bibitem[Rem03]{rempethesis}
Lasse Rempe, \emph{Dynamics of {E}xponential {M}aps}, Ph.D. thesis,
  Christian-Albrechts-Universit\"at zu Kiel, 2003.

\bibitem[RS09]{rempe-schleicher}
Lasse Rempe and Dierk Schleicher, \emph{Bifurcations in the space of
  exponential maps}, Invent. Math. \textbf{175} (2009), no.~1, 103--135.

\end{thebibliography}

\bigskip

\noindent {\bf Weiwei Cui}\\
 Research Center for Mathematics and Interdisciplinary Sciences\\
Frontiers Science Center for Nonlinear Expectations, Ministry of Education\\
Shandong University, Qingdao, 266237, China.

\smallskip

\noindent{weiwei.cui@sdu.edu.cn}

\bigskip

\medskip

\noindent {\bf Jun Wang}\\
School of Mathematical Sciences\\
Fudan University, Shanghai, 200433, China.

\smallskip

\noindent{majwang@fudan.edu.cn}

\end{document}